\documentclass[11pt]{amsart}
\usepackage[utf8]{inputenc}
\usepackage[T1]{fontenc}
\usepackage[colorlinks]{hyperref}
\usepackage{mathtools, amsthm, amsfonts, amssymb}
\usepackage{accents} 
\usepackage{microtype}

\hypersetup{
  linkcolor=[rgb]{0.3,0.3,0.6},
  citecolor=[rgb]{0.2, 0.6, 0.2},
  urlcolor=[rgb]{0.6, 0.2, 0.2}
}

\usepackage{latexsym}
\mathtoolsset{showonlyrefs,showmanualtags}
\usepackage{thmtools}
\usepackage{mathrsfs}
\usepackage{textcomp}
\usepackage{graphicx}
\usepackage{setspace}
\usepackage{nicefrac}
\usepackage{enumerate}
\usepackage{wasysym}
\usepackage[normalem]{ulem}
\usepackage{upgreek}

\usepackage{paralist}
\usepackage{xcolor}
\usepackage{etoolbox}
\usepackage{mathdots}
\usepackage{wrapfig}
\usepackage{floatflt}
\usepackage{tensor} 
\usepackage{lscape}
\usepackage{stmaryrd}
\usepackage[enableskew]{youngtab}
\usepackage{ytableau}
\usepackage{pifont}

\usepackage{tikz}
\usetikzlibrary{arrows}
\usetikzlibrary{matrix}
\usetikzlibrary{positioning}
\usetikzlibrary{calc}

\usepackage[all,cmtip]{xy}
\usepackage[labelfont=bf, font={small}]{caption}[2005/07/16]

\usepackage{xcolor}
\definecolor{red}{rgb}{1,0,0}

\newcommand{\vvirg}{ , \dots , }
\newcommand{\ootimes}{ \otimes \cdots \otimes }

\newcommand{\ttimes}{ \times \cdots \times }

\newcommand{\contract}{\rotatebox[origin=c]{180}{ \reflectbox{$\neg$} }}

\newcommand{\textsum}{{\textstyle \sum}}

\newcommand{\bfP}{\mathbf{P}}

\newcommand{\bfe}{\mathbf{e}}

\newcommand{\bfm}{\mathbf{m}}
\newcommand{\bfn}{\mathbf{n}}

\newcommand{\bfu}{\mathbf{u}}
\newcommand{\bfv}{\mathbf{v}}

\newcommand{\calC}{\mathcal{C}}

\newcommand{\calE}{\mathcal{E}}

\newcommand{\calI}{\mathcal{I}}

\newcommand{\calL}{\mathcal{L}}

\newcommand{\scrS}{\mathscr{S}}

\newcommand{\scrV}{\mathscr{V}}

\newcommand{\bbC}{\mathbb{C}}

\newcommand{\bbN}{\mathbb{N}}

\newcommand{\bbP}{\mathbb{P}}

\newcommand{\bbR}{\mathbb{R}}

\newcommand{\bbX}{\mathbb{X}}

\newcommand{\frakS}{\mathfrak{S}}

\newcommand{\rmQ}{\mathrm{Q}}
\newcommand{\rmR}{\mathrm{R}}

\renewcommand{\phi}{\varphi}
\newcommand{\eps}{\varepsilon}

\renewcommand{\hat}[1]{\widehat{#1}}
\renewcommand{\bar}[1]{\overline{#1}}

\DeclareMathOperator{\codim}{codim}

\DeclareMathOperator{\Hom}{Hom}

\DeclareMathOperator{\trace}{trace}

\DeclareMathOperator{\End}{End}

\DeclareMathOperator{\Sym}{Sym}

\DeclareMathOperator{\Ann}{Ann}

\newcommand{\stabrank}{\mathrm{stabR}}

\newcommand{\Hilb}{\mathrm{Hilb}}

\newcommand{\bfVNP}{\mathbf{VNP}}
\newcommand{\bfVP}{\mathbf{VP}}
\newcommand{\bfNP}{\mathbf{NP}}

\DeclareMathOperator{\dist}{dist}

\newcommand{\uR}{\underline{\rmR}}
\newcommand{\uQ}{\underline{\rmQ}}

\DeclareMathAccent{\wtilde}{\mathord}{largesymbols}{"65}

\newcommand{\cR}{\mathrm{cR}}

\newcommand{\CNOT}{\mathit{CNOT}}

\newcommand{\MaMu}{\mathbf{MaMu}}

\newcommand{\IMM}{\mathrm{IMM}}

\newcommand{\GL}{\mathrm{GL}}

\newcommand{\Mat}{\mathrm{Mat}}

\newcommand{\calStab}{\mathcal{S}\mathit{tab}}

\newcommand{\calTNS}{\mathcal{T\!N\!S}}
\newcommand{\calUTNS}{\mathcal{U\!T\!N\!S}}

\usepackage[top=3cm,bottom=3cm,left=3cm,right=3cm]{geometry}

\title[Geometry of Tensors]{Geometry of Tensors:\\Open problems and research directions}
\author{Fulvio Gesmundo}

\address[Fulvio Gesmundo]{Saarland University, Saarbr\"ucken, Germany}
\email{gesmundo@cs.uni-saarland.de}
\newtheorem{problem}{Problem}

\newcommand{\red}{\mathrm{red}}
\newcommand{\ext}{\mathrm{ext}}
\newcommand{\reg}{\mathrm{reg}}
\newcommand{\EC}{\mathrm{EC}}
\newcommand{\EDdeg}{\mathrm{EDdeg}}
\newcommand{\Ch}{\mathrm{Ch}}
\newcommand{\aEDdeg}{\mathrm{aEDdeg}}
\newcommand{\immC}{\mathrm{immC}}
\newcommand{\uimmC}{\underline{\mathrm{immC}}}
\newcommand{\VWaring}{\mathbf{VW}\!\mathsf{aring}}
\newcommand{\bfVBP}{\mathbf{VBP}}
\newcommand{\VChow}{\mathbf{VC}\!\mathsf{how}}
\newcommand{\uVWaring}{\bar{\mathbf{VW}\!\mathsf{aring}}}
\newcommand{\uVBP}{\bar{\mathbf{VBP}}}
\newcommand{\uVChow}{\bar{\mathbf{VC}\!\mathsf{how}}}

\begin{document}

\begin{abstract}
This is a collection of open problems and research ideas following the presentations and the discussions of the AGATES Kickoff Workshop held at the Institute of Mathematics of the Polish Academy of Sciences (IMPAN) and at the Department of Mathematics of University of Warsaw (MIM UW), September 19-26, 2022.
\end{abstract}

\maketitle

The AGATES Kickoff Workshop\footnote{\url{https://agates.mimuw.edu.pl/index.php/agates/kickoff-workshop}} covered topics ranging from Classical Algebraic Geometry and Representation Theory to Theoretical Computer Science, Quantum Physics and Data Science. Twelve speakers from different fields presented open problems in their area of study, with the goal of providing research directions that the participants could pursue during the rest of the AGATES semester, encouraging new collaborations and the use of tensor geometry in different areas. 

In this document, we organize the research directions into four subjects which reflect those covered during the AGATES Summer School held the week prior to the Kickoff Workshop. In each section, we provide brief information on the background of the corresponding topic. However, this document is by no mean self-contained: we provide some useful references which would help the interested reader to develop a general understanding of the topic, and of the state of the art. It is clear that the references contain much more than what is needed to understand the open problems described below; our suggestion is that the reader should go through these notes, and only later supplement the theory they need to completely understand a problem that piques their interest.

\begin{itemize}
 \item We assume knowledge of algebraic geometry at a basic level. The definition of variety, dimension, degree, the basic properties of the ideal-variety correspondence, and some undergraduate level commutative algebra is assumed. We refer to \cite{Harris:AlgGeo} for the basics of the required projective geometry.
 \item Central in these notes, especially in their first part, is the theory of tensor rank, and some of its variants. \cite{BCCGO:HitchhikerGuide,Lan:TensorBook} cover all the required theory and much more. An introduction for the theory of secant varieties, from a geometric point of view, can be found in \cite{CarGriOed:FourLecturesSecantVarieties,One:NotesAgates}.
 \item We refer to \cite{Ven:NotesAgates} for an introduction to the Spectral Theory of Tensors. The topic is covered extensively from the point of view of tensor analysis in \cite{QiLuo:TensorAnalysis}. More recent developments, together with a study of this topic from a more geometric point of view, can be found in \cite{Sodo:Thesis,Tur:Thesis}.
 \item Geometric and algebraic methods in Theoretical Computer Science and Complexity Theory are covered in detail in \cite{Lan:GeometryComplThBook}. \cite{BlaIk:IntroGCT} covers an introduction to geometric complexity theory, which contains essentially all the material needed to start working on it. \cite{BuClSho:Alg_compl_theory} is a complete textbook on algebraic complexity theory.
 \item The geometry of tensors in quantum physics is extremely diverse. We refer to \cite{Orus:PracticalIntroTensorNetworks,SilTscGerJunJasRizMon:TNAntology} for an introduction to tensor networks. A more geometric point of view is taken in \cite{Stef:Thesis,Sey:MPSgeomIT}. 
\end{itemize}

We point out that the role of tensors is central in geometry and in many related fields and there are a number of related topics which offer many research directions. These notes represent a small cross section. At the end of each section, we propose some additional topics, which have not been explored in detail during the AGATES Kickoff Workshop, but we believe they might be of interest for the readers of these notes as they present numerous open directions and applications.

\subsection*{Acknowledgements} This work is supported by  the Thematic Research Programme ``Tensors: geometry, complexity and quantum entanglement'', University of Warsaw, Excellence Initiative -- Research University and the Simons Foundation Award No. 663281 granted to the Institute of Mathematics of the Polish Academy of Sciences for the years 2021-2023. I would like to thank Weronika Buczy{\'n}ska, Jaros{\l}aw Buczy\'nski, Francesco Galuppi, Joachim Jelisiejew for organizing the AGATES semester and providing a great research environment in Warsaw. I also thank Maria Chiara Brambilla, Mario Kummer, Benjamin Lovitz, Alessandro Oneto, Giorgio Ottaviani, Tim Seynnaeve, Daniel Stilck Fran{\c{c}}a, Vincent Steffan, Nick Vannieuwenhoven, Emanuele Ventura for the inspiring discussions we had during the AGATES semester and their suggestions in the writing of this document.

\subsection*{Notation} 
We mainly work over complex numbers. In some sections, we mention connection to real geometry; we say explicitly when we do so. Given a vector space $V$, $\bbP V$ denotes the projective space of lines in $V$; $V^*$ denotes the dual space of $V$. For spaces $V ,W$, denote by $V \otimes W$ the tensor product of $V$ and $W$; it can be equivalently be identified with the space of linear maps $\Hom(V^*,W) \simeq \Hom(W^*,V)$ or with the space of bilinear map $V \times W \to \bbC$. For $d \geq 0$, let $V^{\otimes d}$ denote the $d$-th tensor power of $V$. Then $S^d V \subseteq V^{\otimes d}$ denotes the subspace of symmetric tensors, which are those tensors invariant under the action of the symmetric group $\frakS_d$ acting by permutation on the tensor factors; $S^d V$ can be identified with the space of homogeneous polynomials of degree $d$ on $V^*$. For a subset $S \subseteq V$, let $\langle S \rangle$ denote the linear span of $S$; if $S \subseteq \bbP V$, the same notation refers to the projective linear span.

\section*{Part 1: Classical Geometry} 
\setcounter{section}{1}
\setcounter{subsection}{0}
Let $V_1 \vvirg V_k$ be vector spaces of dimension $n_1 \vvirg n_k$ respectively and let $d_1 \vvirg d_k$ be non-negative integers. The \emph{Segre-Veronese} embedding of multidegree $(d_1 \vvirg d_k)$ is the map
\[
 v_{d_1 \vvirg d_k} : \bbP V_1 \ttimes \bbP V_k \to \bbP (S^{d_1} V_1 \ootimes S^{d_k} V_k);
 \]
 its image, denoted $\scrV_{d_1 \vvirg d_k}$, is the Segre-Veronese variety of partially symmetric tensors of rank one
 \[
  \scrV_{d_1 \vvirg d_k} = \{ \ell_1^{d_1} \ootimes \ell_k^{ d_k} : \ell_j \in V_j\} \subseteq \bbP (S^{d_1} V_1 \ootimes S^{d_k} V_k).
 \]
When $k =1$, this is the Veronese variety of powers of linear forms. If $(d_1 \vvirg d_k) = (1 \vvirg 1)$, this is the Segre variety of decomposable tensors.

Given a variety $X \subseteq \bbP V$, and an element $p \in \bbP V$, the $X$-rank of $p$ is 
 \begin{equation}\label{eqn: X rank}
  \rmR_X(p) = \min \{ r : p \in \langle x_1 \vvirg x_r \rangle \text{ for some } x_1 \vvirg x_r \in X\};
 \end{equation}
 the $r$-th secant variety of $X$ is 
 \[
  \sigma_r(X) = \bar{\{ p \in \bbP V : \rmR_X(p) \leq r \}}
 \]
where the closure is equivalently in the Zariski or the Euclidean topology. The \emph{border $X$-rank} of $p$ is 
\[
 \uR_X(p) = \min \{ r : r \in \sigma_r(X)\}.
\]
Clearly $\uR_X(p) \leq \rmR_X(p)$ and for most algebraic varieties $X$ there are examples where the inequality is strict.

\subsection{Dimension of secant varieties}
A great deal of research has been spent in the study of secant varieties and of particular interest is the study of their dimension. A straightforward parameter count shows that for a variety $X \subseteq \bbP V$ one has 
\[
 \dim \sigma_r(X) \leq \min \{ r \dim X + r-1, \dim V-1\};
\]
we say that $\sigma_r(X)$ has the \emph{expected dimension} if equality holds; we say that $X$ is $r$-defective, or that $\sigma_r(X)$ is defective, if the inequality is strict.

In general, little is known about dimension of secant varieties. It is a classical fact that curves are never defective \cite{Palatini:SuperficieAlg}; defective surfaces were classified in the the early XX century \cite{Severi,Terra:Superficie}; threefolds with defective second secant variety were classified in \cite{Sco:Determinazione3foldsDef} whereas \cite{ChiCil:ClassifictionDefectiveThreefolds} classifies defective threefolds in general; fourfolds with defective second secant variety have been classified recently \cite{ChiCilRus:DefectiveFourfolds}. In the setting of Segre-Veronese varieties, more is known. Veronese varieties are not defective, except for a finite list of known examples, classified in \cite{AlHir:Poly_interpolation_in_several_variables}. For more general choices of dimension and multi-degree less is known. A summary of the state-of-the-art can be found in \cite{BCCGO:HitchhikerGuide,GalOne:SecantNonDefCollision}. The problem is open in general and it is typically believed that the known defective cases are the only existing ones. 

We propose the following very general research direction:
\begin{problem}\label{prob: defective Segre veronese}
 Classify defective Segre-Veronese varieties and determine their dimension.
\end{problem}
Via Terracini's Lemma \cite[Lemma 1]{BCCGO:HitchhikerGuide}, \autoref{prob: defective Segre veronese} can be rephrased in terms of speciality of certain line bundles on the product of multiprojective spaces. Given a line bundle $\calL$ on a variety $X$ and a zero-dimensional scheme $\bbX \subseteq X$, we say that $\bbX$ imposes independent conditions on $\calL$ if 
\[
 \dim H^0 ( \calI_\bbX \otimes \calL) = \max\{ \dim H^0(\calL) - \deg(\bbX), 0\}.
\]
We say that a line bundle $\calL$ is special if there is $r \geq 1$ such that $r$ double points with generic supports in $X$ do not impose independent conditions on $\calL$. Via Terracini's Lemma, this condition is equivalent to the fact that the $r$-th secant variety of $X$ is defective \cite{BCCGO:HitchhikerGuide}.
\begin{problem}
Classify special line bundles on the product of projective spaces. More generally, classify special line bundles on interesting classes of varieties.
\end{problem}
In general, one may ask whether other $0$-dimensional schemes, different from the union of double points, impose independent condition on certain line bundles. This is, for instance, the object of the Segre-Harbourne-Gimigliano-Hirschowitz Conjecture, which, roughly speaking, predicts that fat points with generic support impose independent conditions except for a classified list of known cases; related conjectures appear in \cite{Nagata:14Hilbert,LafUga:ClassSpecialLinSystemP3}. The case of schemes contained in a union of double points in $\bbP^n$ is studied in \cite{BraOtt:PartialPolyInt}; the case of curvilinear schemes is studied in \cite{CatGim:curvilinearP2}. For a small number of points in $\bbP^n$, the situation is understood, see e.g. \cite{BraDumPos:NotionSpecLinSyst,BraDumPos:EffectiveConeBlowup}.

In principle, one can study the case of schemes with arbitrary supports. A classification in this case is expected to be out of reach. However, specializing the question to the case of $0$-dimensional scheme of degree $2$ makes the problem much more tractable. This case has important connections to the study of Terracini loci \cite{BalChi:TerraciniLocus}, which are varieties associated to the secant varieties and control (at least partially) their singular loci.
\begin{problem}\label{probl: Terracini}
 Given a variety $X$, a line bundle $\calL$ on $X$ and an integer $r$, determine whether there exists a $0$-dimensional scheme $\bbX \subseteq X$ which is union of disjoint schemes of degree $2$, such that the reduced scheme $\bbX^{\red}$ imposes  independent conditions on $\calL$ but $\bbX$ does not.
\end{problem}
For small $r$, several classical results translate to an answer to \autoref{probl: Terracini}. For $r = 2$ the case of curves amounts to characterizing the \emph{edge variety} of the variety $X$; this is well understood for curves, and partial results are known for surfaces  \cite{RanStu:ConvexHullSpaceCurve,RanStu:ConvexHullVar,MerRanSin:ConvexHullSurfaces}. For small $r$, there are several result in the setting of Veronese and Segre-Veronese varieties \cite{BalChi:TerraciniLocus, BalBerSan:TerraciniLocusThreePts,BalVen:AmpleTerracini, ChiGes:DecompTerraciniCubic}. An unrestricted version of \autoref{probl: Terracini} seems hard to approach, but even partial answers for special classes of varieties would be valuable.

\subsection{Generalized additive decompositions and $0$-dimensional schemes}

In order to investigate the relations between rank and border rank of tensors, several related notions of rank have been introduced in the literature. We mention one of them which is important for this section.

Given a variety $X \subseteq \bbP V$ and a point $p \in \bbP V$, the \emph{cactus $X$-rank} of $p$ \cite{BerRan:CactusRankCubicForm,BuczBucz:SecantVarsHighDegVeroneseReembeddingsCataMatAndGorSchemes} is 
\[
 \cR_X(p) = \min\{ r : p \in \langle \bbX \rangle \text{ for some $0$-dimensional scheme $\bbX \subseteq X$, with $\deg(\bbX) = r $}\}.
\]
Note that the definition of cactus rank reduces to the one of rank under the additional hypothesis that the scheme $\bbX$ is smooth, namely a union of $r$ distinct points; in particular $\cR_X(p) \leq \rmR_X(p)$. Of particular interest for this section is the case where $X = v_d(\bbP V) \subseteq \bbP S^d V$ is a Veronese variety. In this case, apolarity theory guarantees that $f \in \langle v_d(\bbX ) \rangle$ if and only the defining ideal $I_\bbX \subseteq \Sym(V^*)$ of $\bbX$ is contained in the \emph{apolar ideal} $\Ann(f)$ of $f$; we refer to \cite[Sec. 2.1.4]{BCCGO:HitchhikerGuide} for an introduction to the theory. We say that $\bbX \subseteq \bbP V$ is apolar to $f$ if $I_\bbX \subseteq \Ann(f)$. We say that $\bbX$ is non-redundant if there is no proper subscheme $\bbX' \subseteq \bbX$ which is apolar to $f$. We say that $\bbX$ is minimal if $\deg(\bbX) = \cR_X(f)$.

In the setting of symmetric tensor, one defines a generalized additive decomposition (GAD) of $f \in S^d V$ to be an expression of the form
\begin{equation}\label{eqn: GAD}
 f = \ell_1^{d-k_1} g_1 + \cdots + \ell_s ^{d-k_s} g_s
\end{equation}
with $\gcd(g_i,\ell_i)=1$ and $\ell_1 \vvirg \ell_s$ distinct. A Waring rank decomposition of $f$ is a special GAD for $f$; similarly, the trivial decomposition $f = f$ is a GAD for $f$, as well. 

One can naturally associate a zero dimensional scheme $\bbX \subseteq \bbP V$ to the GAD of a homogeneous polynomial $f \in S^d V$. The construction is technical and we refer to \cite{BerJelMarRan:PolyGivenHilbertFunction} for the details; its degree gives a notion of \emph{length} for the corresponding GAD. 

Given $f \in S^d V$, the relation between the schemes apolar to $f$ and the schemes associated to GADs for $f$ is not entirely understood. In general, it is true that if $\bbX$ is a scheme associated to a GAD for $f$, then $\bbX$ is apolar to $f$ \cite{BerBraMou:ComparisonDifferentNotionsRanksSymmetricTensors,BerJelMarRan:PolyGivenHilbertFunction}. The converse does not necessarily hold: there are schemes apolar to $f$ which are not induced by a GAD \cite{Tauf:Presentation}. However, a partial converse is true: If $\bbX$ is a scheme apolar to $f \in S^d V$ then there exists a degree $D$ and an \emph{extension} $f^{\ext} \in S^D V$ of $f$ (namely an antiderivative with respect to some differential operator) such that $\bbX$ is associated to a GAD for $f^{\ext}$ \cite{BerBraMou:ComparisonDifferentNotionsRanksSymmetricTensors}. An important parameter that controls these relations is the Castelnuovo-Mumford regularity of the scheme. We refer to \cite{Eisenbud:SyzygyBook} for the basic properties of this notion. If $\bbX$ is a scheme of regularity $\reg(\bbX) \leq d$, then $\bbX$ is apolar to $f$ if and only if it contains a subscheme $\bbX'$ associated to a GAD for $f$. Given $f$, it is however not clear what is the maximum possible $\reg(\bbX)$ for a scheme associated to a GAD for $f$ or for a (minimal) scheme apolar to $f$.

\begin{problem}
Let $f \in S^d V$ and let $\bbX$ be a minimal apolar scheme for $f$. Is $\reg(\bbX) \leq d$? 
\end{problem}
It is known that the bound $\reg(\bbX) \leq d$ does not necessarily hold if $\bbX$ is an irredundant scheme apolar to a form $f \in S^d V$. 
\begin{problem}
Determine the maximum possible value of $\reg(\bbX)$ for the regularity of an irredundant scheme apolar to a form $f \in S^d V$. 
\end{problem}
Moreover, schemes associated to certain particular GADs are regular in degree $d$: this is the case, for instance, of reduced schemes (which are associated to Waring decompositions) and of schemes of degree $2$ (associated to tangential decompositions) \cite{BerTau:WaringTangCactusDecomp}. 
\begin{problem}
For $f \in S^d V$, let $\bbX$ be the scheme associated to a GAD of minimal length. Is $\reg(\bbX) \leq d$? In general, determine what conditions on $f$ guarantee the existence of a GAD whose associated scheme has low regularity. 
\end{problem}

\subsection{Degree of secant varieties}

A classical topic in algebraic geometry is determining the degree of secant varieties. Of particular interest is the case of secant varieties which are hypersurfaces. For instance, in \cite{KumSin:HyeprbolicSecantVars}, it was observed that if $X$ is a real algebraic variety, whose complexification $X_\bbC$ has a secant variety which is a hypersurface, then the defining equation of this hypersurface is a candidate example for a hyperbolic polynomial whose cone cannot be described by a linear matrix inequality. In \cite{ChrGesJen:BorderRankNonMult}, varieties having a secant variety which is a hypersurface were used to produce examples of strict submultiplicativity of rank and border rank under the Segre product.

A lower bound for the degree of the second secant variety is given in \cite{CilRus:VarsMinimalSecant}. If $X \subseteq \bbP V$ is a non-degenerate variety, and $\codim \sigma_2(X) = c$, then $\deg(\sigma_2(X)) \geq \binom{c+2}{2}$. A classical application of the double point formula \cite[Sec. 2.4]{EisHar:3264} allows one to compute the degree of secant varieties of curves: If $C \subseteq \bbP V$ is a non-degenerate curve of degree $d$ and genus $g$, then $\deg(\sigma_2(X)) = \frac{(d-1)(d-2)}{2} - g$. In \cite{CilRus:VarsMinimalSecant}, a more involved use of the double point formula, in connection with the theory of weakly defective varieties, was used to classify varieties $X$ such that $\deg(\sigma_2(X))$ attains the lower bound $\binom{c+2}{2}$.

For higher secant varieties, little is known. Some special cases are entirely or almost entirely understood: we mention rational normal curves and elliptic normal curves \cite{Fis:GenusOneCurvesPfaffians}, and sporadic examples of secant varieties of some Segre-Veronese varieties \cite{Str:RelativeBilComplMatMult,LanMan:IdealSecantVarsSegre,Ottav:InvariantRegardingWaringProblemCubicPolynomials}. 

We propose two problems in this direction:
\begin{problem}
Provide lower bounds for the degree of higher secant varieties, in the spirit of \cite{CilRus:VarsMinimalSecant}.
\end{problem}
\begin{problem}
 Classify varieties of small dimension having second and (or) third secant variety of minimal degree. More generally, given $r$, produce a variety $X$ such that $\sigma_r(X)$ is a hypersurface of the smallest possible degree.
\end{problem}

\subsection{Other topics}

The recent \cite{BuczBucz:ApolarityBorder} developed a \emph{border} version of apolarity theory, based on the multigraded Hilbert scheme introduced in \cite{HaiStu:MultigradedHilbertSchemes}. The study of limits of saturated ideals presents interesting connections to deformation theory \cite{JelMan:LimitsSaturated,Man:IdentiLimitsSaturated}, and offers a number of new interesting open directions. 

Moreover, the research on equations for secant varieties remains an important topic with many open problems. New directions can be explored in connection with the mentioned border apolarity, in the spirit of \cite{ConHarLan:LowerBoundsMatMul33det}.

\section*{Part 2: Spectral Theory of Tensors}
\setcounter{subsection}{0}
\setcounter{section}{2}

Spectral theory of tensors is a topic lying across algebraic geometry and optimization theory. Indeed, it is the natural generalization to the tensor setting of the classical spectral theory of matrices, with its connections to semidefinite programming. The optimization theory is often done over the real numbers, even though the geometry is best done in the complex setting. 

The (real) Euclidean distance plays an important role. Hence, given a complex space $V$, one usually fixes a \emph{real subspace} $V_\bbR$, and a full rank quadratic form (i.e., an inner product) $\langle - , - \rangle \in S^2 V^*$ with the property that its restriction to $V_\bbR$ is positive definite. Then $\langle - , - \rangle$ can be used to define a Euclidean distance $\Vert - \Vert$ on $V_\bbR$, as well as an identification $V \simeq V^*$ of $V$ with its dual space. 

Given $V$ as above, and a tensor $T \in  V^{\otimes d}$, an eigenvector of $T$ is an element $v \in V$ such that $f \underset{2\vvirg d}{\contract}  v ^{\otimes (d-1)} = \lambda v$ for some $\lambda \in \bbC^\times$; here $\underset{2\vvirg d}{\contract}$ denotes the tensor contraction on all but the first tensor factor, after the identification of $V$ with $V$. It is clear that since the tensor $v^{\otimes (d-1)}$ is symmetric, the eigenvectors only depend on the partially symmetric component of $T$ in the space $V \otimes S^{d-1} V$. Similar definitions can be given by choosing the contraction on different tensor factors; this gives rise to the \emph{eigencompatibility variety} discussed in \autoref{sec: eigencompatibility}.

The most immediate generalization from the matrix setting to the tensor setting is for symmetric tensors. If $f \in S^d V$ is identified with a homogeneous polynomial on $V$, then the definition of eigenvectors for $f$ is equivalent to the condition $\nabla f (v) = \lambda v$ for some $\lambda \in \bbC^\times$; here $\nabla f$ denotes the differential of $f$ or equivalently the vector of its partial derivatives for some choice of coordinates.

Denote by $f|_\Sigma$ the restriction of $f$ to a function on the \emph{unit sphere} $\Sigma_V = \{ v \in V : \langle v, v \rangle = 1\}$: it is easy to see that $v$ is an eigenvector for $f$ if and only if it is a critical point of $f|_\Sigma$ with critical value $\lambda = f(v)$. The spectral norm of $f$ is defined to be $\Vert f \Vert_\infty = \max \{ | f(v) | : v \in \Sigma_V\}$ and it coincides with the largest modulus of a critical value of $f$. When $f$ is a quadratic form, or equivalently $T$ is a symmetric matrix, these notions specialize to the analogous notions for matrices.

\subsection{Eigenvectors and eigenschemes of symmetric tensors}

The \emph{eigenscheme} of a homogeneous polynomial $f \in S^d V$ is the subscheme $\bbX_f$ of $\bbP V$ defined by the condition $\nabla f ( v) \wedge v = 0$; if $x_0 \vvirg x_n$ are coordinates on $V$, this is equivalent to the vanishing of the $2 \times 2$ minors of the matrix 
\[
\left( 
\begin{array}{ccc}
x_0 & \cdots  & x_n \\ 
\partial_0 f & \cdots & \partial_n f ,
\end{array}
\right)
\]
where $\partial_j = \frac{\partial}{\partial x_j}$ is the partial derivative in the variable $x_j$. This generalizes the notion of the set of eigenvectors of a symmetric matrix. In \cite{ForSib:ComplexDynHigherDim,CartStur:NumberEigenvalsTensor}, it was observed that if $f$ is generic, then its eigenscheme is $0$-dimensional and it consists of $D_{n,d} = \frac{(d-1)^{n+1} - 1}{d-2}$ distinct points, where $\dim V= n+1$. In fact, \cite{CartStur:NumberEigenvalsTensor} shows a stronger condition: whenever $\bbX_f$ is finite, its degree is $D_{n,d}$.  However, little is known about which configurations of points, or more generally which $0$-dimensional schemes, may arise as eigenschemes of tensors. For every $d$, one defines a \emph{variety of eigenconfigurations}
\[
 \calE_{n,d} = \bar{\{ (v_1 \vvirg v_{D_{n,d}}) \in (\bbP V)^{(\cdot D_{n,d})} : \bbX_f = \{v_1 \vvirg v_{D_{n,d}}\} \text{ for some $f \in S^d V$}\}}.
\]
Here $(\bbP V)^{(\cdot D_{n,d})}$ denotes the symmetrized product $(\bbP V)^{\times D_{n,d}} / \frakS_{D_{n,d}}$. More generally, one can define $\calE_{n,d}$ as a subscheme of the Hilbert scheme $\Hilb_{D_{n,d}}(\bbP V)$ of $D_{n,d}$ points in $\bbP V$.

We propose the following very general problem:
\begin{problem} \label{prob: equations eigenconfiguration}
Determine equations for $\calE_{n,d}$. 
\end{problem}
 In \autoref{prob: equations eigenconfiguration}, ``equations'' is to be intended in a broad sense, as Zariski closed conditions. This problem was investigated in \cite{BeGalVen:EigenschemesTernaryTensors,BeGalVen:EqnsTensorEigen,BeoMR:ConfigurationsEigenpoints}, which provides partial results and essentially a full characterization for $n=2,3$.

One can specialize \autoref{prob: equations eigenconfiguration} to the case of non-reduced schemes. In this case, even before determining equations, one can study what are the possible degrees of the irreducible components of the eigenscheme:
\begin{problem}
 Let $\bbX$ be a non-reduced scheme which arises as eigenscheme of a tensor. What are the possible degrees of the irreducible components of $\bbX$? And what are the schemes that can arise as irreducible components?
\end{problem}
For instance, in the case of matrices, it is known \cite{ASS} that all the schemes arising as irreducible components of an eigenscheme are curvilinear.

\subsection{Eigenschemes compatibility}\label{sec: eigencompatibility}

As mentioned before, eigenvectors of (non-symmetric) tensors in ${V}^{\otimes d}$ can be defined with respect to any choice of $d-1$ factors. The $k$-th eigenscheme of a tensor $T \in {V}^{\otimes d}$ is the eigenscheme defined with respect of the factors different from the $k$-th one; more precisely
\[
 \bbX_T^k = \{ v \in \bbP V : T \underset{\hat{k}}{\contract} v^{\otimes d} = \lambda v \text{ for some $\lambda \in \bbC^{\times}$}\};
\]
here $T \underset{\hat{k}}{\contract} v^{\otimes d}$ denotes the contraction against all but the $k$-th factor. This gives rise to $d$ (usually) distinct eigenschemes, that must respect some compatibility condition. The \emph{eigencompatibility variety}, introduced in \cite{ASS}, records the $d$-tuples of schemes which can be eigenschemes of a tensor $T$:
\[
 \EC_{n,d} = \bar{\{ (\bbX_1 \vvirg \bbX_d) \in (\bbP V)^{\cdot D_{n,d}} \ttimes (\bbP V)^{\cdot D_{n,d}} : \text{ there exists $T \in V^{\otimes d}$ with $\bbX_i = \bbX_T^i$} \} };
\]
and similarly to the variety of eigenconfigurations, one can define the eigencompatibility variety in the product of $d$ copies of $\Hilb_{D_{n,d}}(\bbP V)$. The projection on any of the $d$ factors surjects $\EC_{n,d}$ onto $\calE_{n,d}$. The case $d=2$ of matrices is completely understood: its equations are described in \cite[Prop. 3.1]{ASS}, and either projection onto $(\bbP V)^{\cdot D_{n,d}} = (\bbP V)^{n+1}$ is birational. Besides the case of matrices, little is known and already the case of binary forms hides interesting geometric features. We propose a problem similar to \autoref{prob: equations eigenconfiguration}:
\begin{problem}
Determine dimension and equations for $\EC_{n,d}$. 
\end{problem}

\subsection{Singular vectors of tensors}
Similar to eigenvectors, one can generalize to the tensor setting the notion of singular vector of a matrix. Fix $d$ vector spaces $V_1 \vvirg V_d$, each endowed with an inner product $\langle -,- \rangle_i$ which allow one to identify $V_i$ with $V_i$. Given $T \in V_1 \ttimes V_d$ a singular tuple of $T$ is a $d$-tuple $(v_1 \vvirg v_d) \in V_1 \ttimes V_d$ such that 
\[
 T \underset{\hat{k}}{\contract} (v_1 \ootimes \hat{v_k} \ootimes v_d) = \lambda v_k \text{ for every $k = 1 \vvirg d$}.
\]
The connection to optimization theory is similar to the case of eigenvalues: a singular tuple defines a singular point of $T$ regarded as a function on the product of spheres $\Sigma_{V_1} \ttimes \Sigma_{V_d}$. The maximal singular value defines a norm, that is called the spectral norm of the tensor.

The number of singular $d$-tuples for a generic tensor $T \in V_1 \ootimes V_d$ was determined in \cite{FriOtt:SingularVectorTuples} and given in terms of its generating function. If $n_i +1 = \dim V_i$, then a generic tensor $T$ has $c(\bfn)$ singular tuples where $c(\bfn)$ is the coefficient of $h_1^{n_1} \cdots h_{d}^{n_d}$ in the polynomial
\[
\prod_{k = 1}^d \frac{(\textsum_{i \neq k} h_i)^{n_k+1} - h_k^{n_k+1}}{(\textsum_{i \neq k} h_i) - h_k}.
\]
One can define varieties recording sets of singular tuples, similarly to the case of the variety of eigenconfigurations and the eigencompatibility varieties in the symmetric setting. Problems analogous to \autoref{prob: equations eigenconfiguration} can be studied.

Further, a number of open problems have been proposed regarding the existence of tensors whose singular tuples have particular properties, for instance being defined over the real numbers.

\begin{problem}\label{prob: real singular tuples}
Given real vector spaces $V_1 \vvirg V_d$ with $\dim V_i = n_i+1$, does there exist a tensor $T \in V_1 \ootimes V_d$ such that $T$ has $c(\bfn)$ critical points defined over the real numbers?
\end{problem}
In the symmetric setting, an analog of \autoref{prob: real singular tuples} was posed in \cite{ASS}, asking whether there exist real homogeneous polynomials $f \in S^d V$ admitting $D_{n,d}$ critical points. This problem has affirmative answer: in \cite{Kos:FullyRealEigenconfigurations}, it was shown that there exist harmonic polynomials having fully real eigenconfigurations. Little is known in the non-symmetric setting, besides the case of matrices, where the answer is positive.

\subsection{Best low rank approximations}

An important connection between optimization theory and the spectral theory of tensors lies in the characterization of best rank one approximation. A positive definite quadratic form $\langle -,-\rangle$ on a real vector space $V_\bbR$ defines a Euclidean norm $\Vert v \Vert^2 = \langle v,v\rangle$, hence a distance function $\dist(v,w) = \Vert v - w\Vert^2$. 

Given an algebraic variety $X \subseteq V$, and an element $v \in V$, the distance function naturally induces a regular function on $X$ defined by $\dist_X(-,v) : X \to \bbC$. If $v$ is generic, the set of critical points of $\dist_X(-,v)$ on $X$ is finite, it consists of distinct points, and its cardinality does not depend on the choice of $v$; the \emph{Euclidean distance degree} of $X$, denoted $\EDdeg(X)$ is the number of critical points of $\dist(-,v)$ for a generic choice of $v \in V$ \cite{DraHorOttStuTho:EDdegree}. If $X$ coincides with the Zariski closure of the set of its real points, one can show that at least one of these critical points is real: the real point $\xi$ realizing the lowest possible value of $\dist(v,\xi)$ for $\xi \in X$ is the best approximation of $v$ on $X$.

In the setting of symmetric tensors, the inner product $\langle -,-\rangle$ on a vector space $V$ of dimension $n+1$ naturally induces an inner product on $S^d V$, often called the Frobenius inner product. Denote by $\scrV_{d,n}$ the affine cone over the Veronese variety $v_d(\bbP V) \subseteq \bbP S^d V$. It is easy to prove that for $f \in S^d V$, a scalar multiple of the element $v^{d} \in \scrV_{n,d}$ is a critical point of $\dist_{\scrV_{n,d}}(-,f)$ if and only if $[v]$ is an eigenvector of $f$. In particular, the results of \cite{CartStur:NumberEigenvalsTensor} guarantee that $D_{n,d} = \EDdeg(\scrV_{n,d})$. 

Analogously, in the case of arbitrary tensors, inner products $\langle - , - \rangle$ on spaces $V_i$ with $\dim V_i = n_i+1$ induce naturally a Frobenius inner product on $V_1 \ootimes V_d$. Denoting by $\scrS_{n_1 \vvirg n_d}$ the affine cone over the Segre variety of rank one tensors in $\bbP ( V_1 \ootimes V_d)$, one has a natural correspondence between singular tuples of a tensor $T$ and critical vectors of the distance function $\dist_{\scrS_{n_1 \vvirg n_d}} ( - , T)$. In this case, the number $c(\bfn)$ from \cite{FriOtt:SingularVectorTuples} coincides with the Euclidean distance degree $\EDdeg(\scrS_{n_1 \vvirg n_d})$.

Notice that the above characterizations of the Euclidean distance degree in the tensor setting concern the distance function from the affine cone over the projective variety of rank one tensors. One can define an affine version of the Euclidean distance degree as follows: Let $\hat{X}$ be the cone over a projective variety $X \subseteq \bbP V$ and let $H \subseteq V$ be a generic affine hyperplane; let $X^H = \hat{X} \cap H$, which is an affine variety in $H$. The Euclidean distance degree of $X^H$ does not depend on the choice of $H$, and is called the affine Euclidean distance degree of $X$, denoted $\aEDdeg(X)$. In general, $\EDdeg(\hat{X})$ and $\aEDdeg(X)$ are different. In the tensor setting, this notion has been studied in \cite{BoePetStu:MarginalIndependenceModels} in the context of probability tensors, that are tensors $T \in V_1 \ootimes V_d$ lying in the affine hyperplane $\sum_{(i_1 \vvirg i_d) \in [n_1] \ttimes [n_d]} t_{i_1 \vvirg i_d} = 1$. It is conjectured \cite[Conj. 20]{BoePetStu:MarginalIndependenceModels} that if $T$ is a probability tensor with nonnegative coordinates, then, among its $\aEDdeg (\scrS_{n_1 \vvirg n_d})$ critical points, there is exactly one with real entries. Hence, we pose the following problem
\begin{problem}
Determine the possible numbers of real elements among the $\aEDdeg(\scrS_{n_1 \vvirg n_d})$ critical points of a probability tensor.
\end{problem}

We point out that the problem of determining the Euclidean distance degree for varieties different from rational homogeneous varieties is essentially open, other than for special cases. For instance, it is completely open for secant varieties of Segre and Veronese varieties, already in small dimension. We propose a problem in a restricted setting:
\begin{problem}
 Determine $\EDdeg(\sigma_2(v_d(\bbP^1)))$, the Euclidean distance degree for the secant variety of the rational normal curve in $\bbP ^d$. 
\end{problem}

\subsection{Other topics}
The spectral norm of a tensor is a quantity of interest in optimization theory and, unlike the case of matrices, very few methods to compute it are known. In particular, an important quantity of interest is the \emph{best rank one approximation ratio} between the Frobenius norm and the spectral norm. This has been studied from a geometric point of view in \cite{AgrKozUsc:ChebyshevPoly} and it is known in very few cases. 

This topic can also be studied in a restricted setting. For instance, in \cite{EisUsc:MaximumRelDistanceBetweenRankTwo}, the distance function restricted to the variety of rank $2$ tensors was studied, and the maximum distance of a (border) rank $2$ tensor from the variety of rank one tensors was determined. 

\section*{Part 3: Geometric methods in Algebraic Complexity Theory}
\setcounter{subsection}{0}
\setcounter{section}{3}

Algebraic tasks such as evaluating polynomials, multiplying matrices or solving linear and polynomial systems are central in several areas of mathematics. Algebraic complexity theory studies the complexity of performing these tasks, employing tools from combinatorics, algebra, probability and geometry. We will focus on the task of evaluating polynomials, typically modeled using \emph{algebraic circuits}: these are directed graphs encoding an algorithm that evaluates a given polynomial. Roughly speaking, a polynomial (more precisely a family of polynomials) is considered easy to evaluate if it admits a small circuit; the meaning of ``small'' often depends on the particular complexity class that one examines. Usually, one restricts to special circuits, reflecting the computational model of interest; in this setting, the fact that a certain family of polynomial belongs to a given complexity class can be translated into membership into some algebraic variety, and it can be studied using tools from geometry and representation theory.

For simplicity, we restrict to the homogeneous setting. A p-family is a sequence $(f_n)_{n \in \bbN}$ where $f_n \in S^{d_n} \bbC^{N_n}$ is a homogeneous polynomial of degree $d_n$ in $N_n$ variables, where $d_n,N_n$ are polynomially bounded as functions of $n$. A complexity class is usually defined to be the set of p-families for which a given measure of complexity is polynomially bounded as functions of the parameter $n$ indexing the sequence.

\subsection{Relations between different measures of complexity}

Recall the definition of $X$-rank with respect to an algebraic variety $X$ from \eqref{eqn: X rank}. Two varieties of interest in complexity theory are the the Veronese variety $ \scrV_{N,d}$, already mentioned before, and the Chow variety of products of linear forms 
\[
 \Ch_{N,d} = \{ \ell_1 \cdots \ell_d \in S^d \bbC^N : \text{ for some } \ell_j \in S^1 V\}.
\]
The Chow rank of a homogeneous polynomial $f \in S^d \bbC^N$, denoted $\rmR_{Ch}(f)$, is its rank with respect to the Chow variety. Correspondingly, one defines complexity classes 
\begin{align*}
 \VWaring &= \{ (f_n) \text{ p-family}: \rmR_{\scrV}(f_n) \in O(n^s) \text{ for some $s$}\};\\
 \VChow &= \{ (f_n) \text{ p-family}: \rmR_{\Ch}(f_n) \in O(n^s) \text{ for some $s$}\};
\end{align*}
these are, respectively, the classes of sequences of polynomials with polynomially bounded Waring and Chow rank. In the computer science community these are sometimes denoted $\Sigma\Lambda\Sigma$ and $\Sigma\Pi\Sigma$, although these classes do not necessarily require the polynomials to be homogeneous.

Another interesting measure of complexity is the \emph{iterated matrix multiplication complexity} of a polynomial $f \in S^d \bbC^N$, defined as follows: 
\[
 \immC(f) = \min \left\{ r : f = M_1 \cdots M_d , \text{ for some } \begin{array}{l}
                                                             M_1 \in \bbC^N \otimes \Mat_{1 \times r} \\
    M_j \in \bbC^N \otimes \Mat_{r \times r}  \text{ for $j=2 \vvirg d-1$}\\
    M_d \in \bbC^N \otimes \Mat_{r \times 1} \end{array}\right\};
\]
in other words $\immC(f)$ is the smallest $r$ for which $f$ can be written as a product of matrices of size $r$ whose entries are linear forms. In the complexity theory community, $\immC(f)$ is called the \emph{algebraic branching program width} of $f$ and the expression of $f$ as product of matrices is called an algebraic branching program (ABP) for $f$. The complexity class $\bfVBP$ is 
\[
 \bfVBP = \{ (f_n) \text{ p-family} : \immC(f_n) \in O( n^s) \text{ for some $s$}\}.
\]
For a homogeneous polynomial $f$, it is clear that $\rmR_{\scrV}(f) \geq \rmR_{\Ch}(f) \geq \immC(f)$. Hence one has the inclusions $\VWaring \subseteq \VChow \subseteq \bfVBP$. It is known that these inclusions are strict: the p-family $m_n = x_1 \cdots x_n$ satisfies $\rmR_{\scrV}(m_n) = 2^{n-1}$ \cite{RanSch:RankSymmetricForm}, whereas $\rmR_{\Ch}(m_n) = 1$, showing $(m_n) \in \VChow \setminus \VWaring$; the p-family $\det_n = \det(X_n)$, where $X_n$ is a matrix of variables, satisfies $\rmR_{\Ch}(\det_n) \geq c_n$ where $c_n$ is a super-polynomial function of $n$\cite{GKKS:ArithmeticCircuitsChasmDepthThree,LimSriTav:SuperpolyDepth3} and $\immC(\det_n) \leq O(n^3)$ \cite{Val:Complete_classes_in_algebra,Toda:Classes_of_arithm_circuits_for_det}, showing $(\det_n) \in \bfVBP \setminus \VChow$. 

One can analogously define \emph{border} classes, by replacing the complexity measures $\rmR_{\scrV},\rmR_{\Ch} , \immC$ with their semicontinuous analogue $\uR_{\scrV},\uR_{\Ch} , \uimmC$. We have seen already the definition of border rank. The \emph{border} iterated matrix multiplication complexity is defined by 
\[
 \uimmC(f) = \min\{ r : f = \lim_{\eps \to 0} f_\eps \text{ for some sequence $f_\eps$ such that $\immC(f_\eps) \leq r$}\}.
\]
The corresponding complexity classes are denoted $\uVWaring, \uVChow,\uVBP$. Again, one immediately has $\uVWaring \subseteq \uVChow \subseteq \uVBP$. The same examples as above guarantee that these inclusions are strict. Little is known about the difference between a complexity class and its border analogue.

\begin{problem}
 Let $\calC$ be a complexity class among $\VWaring, \VChow$ or $\bfVBP$. Determine whether $\calC = \bar{\calC}$. 
\end{problem}

In fact, even weaker \emph{debordering} results are of great interest. For instance, the recent \cite{DutDwiSax:DemystifyingBorderDepth3} showed the inclusion $\uVChow \subseteq \bfVBP$; \cite{DGIJL:BorderComplexity} gives an upper bound for $\rmR_{\scrV}(f)$ of the form $\deg(f) \cdot \exp( \uR_{\scrV}(f) ) $. We propose a major open problem in this area:
\begin{problem}
 Determine whether $\uVWaring \subseteq \VChow$. More modestly, determine non-trivial upper bounds for $\uR_{\scrV}(f)$ in terms of $\rmR_{\Ch}(f)$.
\end{problem}

\subsection{Variety of ABPs of small width}\label{subsec: ABP}

Similarly to the setting of secant varieties, the set 
\[
 \IMM_{r,V} = \{ f \in \bbP S^d V: \uimmC(f) \leq r\}
\]
is closed in the Zariski topology, so that it can be regarded as a subvariety of $\bbP S^d V$. This is the variety of homogeneous polynomials admitting an algebraic branching program (ABP) of width at most $r$. When $r = 1$, this coincides with the Chow variety $\Ch_{d,n}$. \cite{GGIL:DegreeRestrictedStrengthDecompsABP} provides some equation for $\IMM_{r,V}$ based on intersection theoretic properties of the hypersurface $\{ f = 0\} \subseteq \bbP V$, in connection with other complexity measures on the space of polynomials, such as slice rank and strength \cite{BikDraEgg:PolyTensBoundedStrength}. \cite{AllWang:PowerABP2} gives a characterization of $\IMM_{2,V}$, based on certain degeneracy conditions of restrictions of $f$. The value of $\dim \IMM_{r,V}$ can be determined using the results of \cite{Ges:Geometry_of_IMM}. We propose to extend these results for a wider range of values of $r$ and $\dim V$.
\begin{problem}
 Determine equations for $\IMM_{r,V} \subseteq \bbP S^d V$ for small values of $r$ and $\dim V$. 
\end{problem}
The geometric complexity theory program (GCT) \cite{GCT1,GCT2} is a proposed approach to studying the relations between (border) complexity classes via representation theoretic \emph{obstructions}. We refer to \cite{BlaIk:IntroGCT} for an extensive introduction to the method which we summarize here briefly. Usually, varieties $X_r$ controlling complexity measures, such as $\sigma_r(v_d(\bbP V))$, $\sigma_r(\Ch_{d,V})$ and $ \IMM_{r,V}$, are $\GL(V)$-varieties, that is varieties closed under the natural action of the general linear group $\GL(V)$. This makes their ideals $\calI(X_r) \subseteq \Sym(S^d V^*)$ and their coordinate rings $\bbC[X_r] = \Sym(S^d V^*) / \calI(X_r)$ into representations for the group $\GL(V)$. In general, if two projective varieties $X,Y$ satisfy $X \subseteq Y$, then there is a surjection $\bbC[Y] \to \bbC[X]$; moreover this surjection is graded, and if $X,Y$ are $\GL(V)$-varieties, it is $\GL(V)$-equivariant. For this reason, given two varieties $X_r$ and $Y_s$ controlling two $\GL(V)$-invariant (border) complexity measures, one can prove $X_r \not \subseteq Y_s$ by showing that there is no surjection $\bbC[Y_s]_\delta \to \bbC[X_r]_\delta$ for some degree $\delta$. Because of the $\GL(V)$-equivariancy, a sufficient condition to guarantee such a surjection cannot exist is the existence of an irreducible $\GL(V)$-representation occurring in $\bbC[X_r]$ with higher multiplicity than in $\bbC[Y_s]$. An irreducible with the property above is called a \emph{multiplicity obstruction} and its existence guarantees $X_r \not \subseteq Y_s$. The geometric complexity theory program proposed this approach, and some of its variants, as a method toward proving $\bar{\bfVNP}\not \subseteq \bar{\bfVP}$, but the same framework can be adapted to separation of other border complexity classes.

However, the results of \cite{IkPa:Rectangular_Kron_in_GCT,BuIkPa:no_occurrence_obstructions_in_GCT,GesIkPa:GCTMatrixPowering} pose strong restrictions on the complexity models and the varieties to which this method can possibly be applied successfully. One setting in which this method has not been yet explored is the one of $\IMM_{r,V}$. 
\begin{problem}
 Implement the GCT approach to the variety of small algebraic branching programs $\IMM_{r,V}$.
\end{problem}

\subsection{Other topics}

In this section, we did not cover any problem regarding the matrix multiplication complexity, a central subject in algebraic complexity theory. The research for lower bounds on the tensor rank of the matrix multiplication tensors has been a major motivation for the research on equations for secant varieties \cite{LanOtt:EqnsSecantVarsVeroneseandOthers,LanOtt:NewLowerBoundsBorderRankMatMult,LanMic:logLowerBoundBorderRank}. The use of geometric methods for upper bounds appears already in \cite{BiCaLoRo:O277ComplexityApproximateMatMult}, which first motivated a systematic study of border rank; the most recent upper bounds \cite{CopperWinog:MatrixMultiplicationArithmeticProgressions,AlmWil:RefinedLaserMethodFMM} are based on what is known as Strassen's laser method \cite{Str:RelativeBilComplMatMult}, a mix of geometric, probabilistic, and combinatorial techniques. Several barriers have been proved for this method, at least in its original form \cite{AlmWil:LimitsAllKnownApproachesMaMu,AmbFilLeG:FastMaMuLimitsCopWinMethod,ChrVraZui:BarriersIrreversibility,BlaLys:SliceRankBlockIrreversStructureTensors}; \cite{HomJelMicSey:BoundsComplAwayFromCW} proposes a modified approach with the goal of circumventing such barriers.

\section*{Part 4: Geometric methods in quantum physics}
\setcounter{subsection}{0}
\setcounter{section}{4}

This final section concerns problems on the geometry of tensors which originate from the study of entanglement in quantum physics and quantum information theory. The problems described in this section, however, have strong connections to other topics as well, such as invariant theory and additive combinatorics; we will briefly mention these connections when relevant.

Briefly, quantum physics models the state of a quantum system of interest as an element of a complex vector space, endowed with a Hermitian inner product. The state of a composite system is described by an element of a tensor space, where each factor is the space corresponding to one of the components of the composite system. The evolution of the system is described by the action of a unitary matrix. We refer to \cite{NielChu:QuantumCompQuantumInf} for an extensive account of the quantum formalism and an introduction to quantum information theory.

From a geometric point of view, we focus on quantum SLOCC operations, which ultimately correspond to transformations on a tensor space $V_1 \ootimes V_d$ induced by the action of $\End(V_1) \ttimes \End(V_d)$. Hence, much of the theory of entanglement can be studied geometrically, in terms of orbits and orbit-closures in tensor spaces. To this end, we introduce the notions of restriction and degeneration of tensors. Given $T,S \in V_1 \ootimes V_d$, we say that $T$ restricts (resp. degenerates) to $S$ if $S \in \End(V_1) \ttimes \End(V_d) \cdot T$ (resp. $S \in \bar{\End(V_1) \ttimes \End(V_d) \cdot T}$). It is a classical fact that the variety $\bar{\End(V_1) \ttimes \End(V_d) \cdot T}$ of the degenerations of $T$ coincides with the orbit-closure of $T$ under the action of $\GL(V_1) \ttimes \GL(V_d)$. We will say that a tensor $T$ restricts or degenerates to a tensor $S$ even if they do not belong to the same tensor space: this  means that the respective spaces are reembedded into a common space where the image of $T$ restricts or degenerates to the image of $S$. For instance, consider the unit tensor of rank $r$ on $d$ factors, defined by $\bfu_d(r) = \sum_{i=1}^r e_i \ootimes e_i \in \bbC^r \ootimes \bbC^r$; the tensor rank $\rmR(T)$ (resp. the tensor border rank $\uR(T)$) of a tensor $T \in V_1 \ootimes V_d$ is characterized as the smallest $r$ such that $\bfu_d(r)$ restricts (resp. degenerates) to $T$.

\subsection{Tensor networks and matrix product states}

Tensor networks define a class of tensors arising via specific tensor contractions which are encoded in the combinatorics of a weighted graph. They are used in quantum many-body physics because they model desirable entanglement properties, which find application in holography, quantum chemistry and machine learning. Specifically, let $\Gamma = (\bfv(\Gamma), \bfe(\Gamma))$ be a graph with set of vertices $\bfv(\Gamma) = \{ 1 \vvirg d\}$ and set of edges $\bfe(\Gamma)$, and let $\bfm = (m_e : e \in \bfe(\Gamma))$ be an assignment of integer weights on the edges of $\Gamma$, called \emph{bond dimensions}. One defines a specific \emph{graph tensor} $T(\Gamma, \bfm)$ in a tensor space $W_1 \ootimes W_d$; the dimension of the factor $W_j$ is determined by the weights of the edges incident to $j$; more precisely $W_j = \bigotimes_{e \ni j} \bbC^{m_e}$. Further, given a set of integer weights $\bfn = (n_j : j = 1 \vvirg n)$ on the vertices of $\Gamma$, called \emph{local dimensions}, one defines a set 
\[
{ \calTNS^\Gamma_{\bfm,\bfn} }^\circ = \{ T \in \bbC^{n_1} \ootimes \bbC^{n_d} : T \text{ is a restriction of $T(\Gamma,\bfm)$}\}
\]
The tensor network variety is $\calTNS^\Gamma_{\bfm,\bfn} = \bar{{\calTNS^\Gamma_{\bfm,\bfn} }^\circ}$ where the closure can be taken equivalently in the Zariski or the Euclidean topology. We refer to \cite{BerDLaGes:DimensionTNS} for the details of this construction and to \cite{ChrVraZui:AsyRankGraph,ChrGesMicZui:BorRankNonAddHigher} for generalizations to the case of hypergraphs.

In \cite{BerDLaGes:DimensionTNS}, an upper bound to the dimension of the variety $\calTNS^\Gamma_{\bfm,\bfn} $ was given and it was shown that such upper bound is sharp in a wide range of bond and local dimensions. Moreover, defining equations for the variety $\calTNS^\Gamma_{\bfm,\bfn} $ were given for some particular cases. We propose the following problem:
\begin{problem}\label{prob: equations TNS}
 Determine equations for the tensor network variety $\calTNS^\Gamma_{\bfm,\bfn} $. Restrict, if necessary, to particular classes of graphs $\Gamma$, or particular assignments of bond and local dimensions.
\end{problem}
Interestingly, in the case where $\Gamma$ is a cycle graph $C_d$, with bond dimensions $m_1 \vvirg m_d$, the graph tensor coincides with the iterated matrix multiplication tensor, which encodes the multilinear map 
\begin{align*}
\MaMu_{\bfm} :  \Mat_{m_1 \times m_2}(\bbC) \ttimes  \Mat_{m_{d-1} \times m_{d}}(\bbC) &\to \Mat_{m_{1} \times m_{d}}(\bbC) \\
(A_1 \vvirg A_{d-1}) &\mapsto A_1 \cdots A_{d-1}.
\end{align*}
In this case, the tensor network variety is known as the variety of \emph{matrix product states} of bond dimension $\bfm$ \cite{PerVerWolCir:MPSrepresentations}, and it is the subject of a rich line of research. It is a non-commutative analog of the variety of algebraic branching programs of small width defined in \autoref{subsec: ABP}; similarly to the case of ABPs, one easily determines a linear-algebraic description of the variety of matrix product states
\[
 \calTNS^{C_d}_{\bfm,\bfn} = \bar{\left\{ \sum_{\substack{(i_1 \vvirg i_d) \in [n_1] \ttimes [n_d]}} \trace(M^{(1)}_{i_1} \cdots M^{(d)}_{i_d}) \cdot e_{i_1} \ootimes e_{i_d} : M^{(k)}_{i_k} \in \Mat_{m_k \times m_{k+1}}(\bbC) \right\}}.
\]
In this case, and more generally in cases where the underlying graph of a tensor network presents particular symmetries, one is interested in a restricted version of the matrix product states: if $m_1 = \cdots = m_d$ and $n_1 = \cdots = n_d$, one considers the subvariety of $\calTNS^{C_d}_{\bfm,\bfn}$ arising from collections of matrices which do not depend on the index $k$. More precisely
\[
 \calUTNS^{C_d}_{m,n} = \bar{\left\{ \sum_{\substack{(i_1 \vvirg i_d) \in [n] \ttimes [n]}} \trace(M_{i_1} \cdots M_{i_d}) \cdot e_{i_1} \ootimes e_{i_d} : M_{i_k} \in \Mat_{m \times m}(\bbC) \right\}} 
\]
This is the variety of \emph{uniform} or \emph{translation invariant} matrix product states with periodic boundary conditions, and their importance in quantum physics is explained in \cite{PerVerWolCir:MPSrepresentations,VWPC:CriticalityAreaLaw}. In \cite{GesLanWal:MPSandQuantumMaxFlowMinCut}, some equations for this variety are obtained in terms of certain rank conditions. More equations are obtained in \cite{CzMiSe:UniformMPSfromAG}, and \cite{DeMotSey:LinearSpanMPS} provides the complete characterization of the space of linear equations in a restricted setting. A central role in these results is played by the \emph{trace algebra}, which allows one to draw strong connections between the geometry of $\calUTNS^{C_d}_{m,n}$ and classical results in invariant theory. Many problems remain open in this area; we propose the following general problem and we refer to \cite{Sey:MPSgeomIT} for further details.
\begin{problem}
Study the geometry of the variety of uniform matrix product states. In particular, use the properties of the trace algebra to systematically determine families of equations.
\end{problem}

\subsection{Tensor subrank}
The subrank of a tensor is a measure of how much a tensor can be diagonalized by taking linear combinations of its slices in all directions; this notion was introduced in \cite{Str:RelativeBilComplMatMult} in the study of matrix multiplication complexity, and it is in a sense dual to the classical notion of tensor rank defined before. Precisely, given $T \in V_1 \ootimes V_d$, the subrank of $T$ is 
\[
 \rmQ(T) = \max \{ q : T \text{ restricts to } \bfu_d(q)\};
\]
Unsurprisingly, the subrank is not semicontinuous, which motivates the definition of the border subrank of $T$, that is 
\[
 \uQ(T) = \max \{ q : T \text{ degenerates to } \bfu_d(q)\}.
\]
Clearly $\rmQ(T) \leq \uQ(T)$ and, by definition, the border subrank is semicontinuous \emph{along degenerations}, in the sense that if $S$ is a degeneration of $T$ then $\uQ(S) \leq \uQ(T)$. However, it is not necessarily semicontinuous along arbitrary curves: for instance, if $Z$ is a generic rank one tensor in $\bbC^q \ootimes \bbC^q$, then $T_\eps = \bfu_d(q) + \eps Z$ satisfies $\uQ(T_\eps) \leq q-1$ for generic $\eps$, but $T_0 = \bfu_d(q)$ so $\uQ(T_0) = q$.

The subrank is a measure of entanglement in quantum information theory. It controls the amount of ``GHZ-type'' entanglement that can be \emph{distilled} from a quantum state. We refer to \cite{VrChr:EntDistGHZShares} for the details.

Usually, lower bounds on the subrank and border subrank are obtained via explicit reductions. On the other hand, we have only few methods for upper bounds, and they usually rely on determining some other notion of rank which bounds from above the subrank and the border subrank; this is the case for geometric rank \cite{KopMosZui:GeomRankSubrankMaMu}, slice rank \cite{SawTao:NoteSliceRank}, tensor compressibility values \cite{LanMic:logLowerBoundBorderRank}, $G$-stable rank \cite{Derk:GStableRank} and a number of other notions. Although these bounds have been useful to determine exactly the border subrank (and the subrank) in many cases of interest, there is evidence that we cannot expect them to be sharp in general. So we propose the following problem.
\begin{problem}
 Determine a genuine method for upper bounds on the subrank and border subrank of a tensor that does not rely on other notions of rank.
\end{problem}

In \cite{DerkMakZui:GenericSubrank}, the value of $\rmQ(T)$ for the subrank for a generic tensor $T \in \bbC^n \ootimes \bbC^n$ was determined, within a small margin of error. For tensors on three factors, it was shown that $\lfloor \sqrt{3n-2} \rfloor - 5 \leq \rmQ(T) \leq \lfloor \sqrt{3n-2} \rfloor$. The value of the border subrank of generic tensors is not known; if $T \in \bbC^n \ootimes \bbC^n$ is generic, it is known $\uQ(T) \leq n-1$, see e.g. \cite{Chang:MaximalBorderSubrank}.
\begin{problem}
 Determine the border subrank of a generic tensor and give conditions for a tensor to have border subrank higher than the generic one.
\end{problem}

\subsection{Stabilizer rank}

The notion of stabilizer rank appears in entanglement theory. It is a measure of the complexity of \emph{preparing} a given quantum state using certain specific quantum operations. Let $X,Y,Z$ be the three Pauli matrices, defined by 
\begin{align*}
 X = \left( \begin{array}{cc} 0 & 1 \\ 1 & 0 \end{array}\right), \quad 
 Y = \left( \begin{array}{cc} 0 & -i \\ i & 0 \end{array}\right), \quad 
 Z = \left( \begin{array}{cc} 1 & 0 \\ 0 & -1 \end{array}\right) .
  \end{align*}
The Pauli group is the subgroup $P_1 \subseteq \GL_2(\bbC)$ generated by $\{ X,Y,Z\}$; it is a group of $16$ elements. The $d$-th Pauli group is the group $P_d \subseteq \GL(\bbC^2 \ootimes \bbC^2)$ image of $P_1 \ttimes P_1 \subseteq \GL_2 \ttimes \GL_2$ via the natural action of $\GL_2 ^{\times d}$ on ${\bbC^2}^{\otimes d}$.

The Clifford group is the normalizer of $P_d$ in the group $U({\bbC^2}^{\otimes d})$ of unitary matrices. It is infinite, but the quotient over its center is finite; denote this quotient by $C_d$. There is a natural action of $C_d$ on $\bbP (\bbC^2 \ootimes \bbC^2)$. To describe a set of generators for $C_d$, we introduce two elements of $\GL_2$
\begin{align*}
H = \frac{1}{\sqrt{2}}\left( \begin{array}{cc} 
            1 & 1 \\ 1 & -1 
           \end{array}\right),  \quad 
 S = \left( \begin{array}{cc} 1 & 0 \\ 0 & i \end{array}\right), \\
  \end{align*}
and one of $\GL(\bbC^2 \otimes \bbC^2)$
\begin{align*}
 \CNOT : \bbC^2 \otimes \bbC^2 &\to \bbC^2 \otimes \bbC^2 \\
 e_0 \otimes e_0 &\mapsto e_0 \otimes e_0 \\
 e_0 \otimes e_1 &\mapsto e_0 \otimes e_1 \\
 e_1 \otimes e_0 &\mapsto e_1 \otimes e_1 \\
 e_1 \otimes e_1 &\mapsto e_1 \otimes e_0 .
\end{align*}
Then $C_d$ is generated by (the classes modulo the center of $U({\bbC^2}^{\otimes d})$ of) the elements $H , S, \CNOT$ acting on any factor, or pair of factors of $\bbC^2 \ootimes \bbC^2$. 

The set of \emph{stabilizer states} is
\[
\calStab_d = \{ U \cdot (e_0^{\otimes d})  \in \bbC^2 \ootimes \bbC^2 : U \in C_d\};
\]
the set $\calStab_d \subseteq {\bbC^2} ^{\otimes d}$ is finite, and it has $2^d \prod_{k=1}^d (2^k +1)$ elements \cite{Gross:HudsonThm} and its linear span is the entire space $\bbC^2 \ootimes \bbC^2$.

The stabilizer rank of a state $T \in \bbC^2 \ootimes \bbC^2$ is defined by 
\[
 \stabrank(T) = \min\{r : T = \alpha_1 T_1 + \cdots +  \alpha_r T_r \text{ for some $T_i \in \calStab_d$, and $\alpha_i \in \bbC$}\};
\]
in other words $\stabrank(T) = \rmR_{\calStab_d}(T)$.

It is in general a hard problem to exhibit explicit states with large stabilizer rank. If $v = e_0 + \sqrt{2} e_1$, it is shown in \cite{PSV:LowerBoundsStabRank} that $\stabrank(v ^{\otimes d}) \geq c d$ for some constant $c$; the slightly weaker bound $\stabrank(v ^{\otimes d}) \geq \frac{d+1}{4\log(d+1)}$ is shown in \cite{LovSte:NewTechStabRank}.
\begin{problem}
 Prove a super linear lower bound for stabilizer rank.
\end{problem}
Since $\calStab_d$ is finite, the set of elements of stabilizer rank bounded by $r$ is an arrangement of $r$-dimensional linear spaces. In \cite{LovSte:NewTechStabRank}, these sets are studied from a geometric and combinatorial point of view. In particular, the \emph{generic stabilizer rank} $\stabrank_d$ is defined to be the smallest $r$ such that $S^d \bbC^2$ is contained in the span of $r$ elements of $\calStab_d$. It is known \cite{QasPasGos:ImprovedUpperBound,LovSte:NewTechStabRank} that $\stabrank_d \leq O(2^{d/2})$ and clearly $\stabrank_d \geq d+1$; moreover, $\stabrank_d$ is known to be superpolynomial as a function of $d$ unless $\bfP = \bfNP$.
\begin{problem}
 Determine non-trivial lower bounds on $\stabrank_d$.
\end{problem}

\subsection{Other topics} Tensor network representations of quantum states of interest are desirable because they allow one to perform tensor operations \emph{locally}, only considering few factors at the time. Throughout the years, several methods to do this efficiently have been proposed. One recent one is based on Metropolis Monte-Carlo contraction, proposed in \cite{ScuWolVerCir:SimulationQuantum,aragones2021classical} in the setting of matrix product states. One can explore whether this method can be extended to other tensor networks, see e.g. \cite{CapRouSti:ModifiedLogSobolev}. 

Another topic of interest, only tangentially touched in these notes, is \emph{Strassen's spectral theory}, which has connections to complexity theory, additive combinatorics and quantum information. There are several open problem in this area, and we refer to \cite{WIgZui:AsymptoticSpectra} for a complete reference.

\bibliographystyle{alphaurl}
\bibliography{tensors.bib}

\end{document}